\documentclass[10pt,a4paper]{amsart}
\usepackage{amssymb,amscd,amsmath,amsthm, young,mathrsfs}%,mathbbol}
\usepackage{mathdots} 
\usepackage[mathcal]{eucal} 
\usepackage{float}
\usepackage[usenames]{color}
\usepackage{soul}
\usepackage{booktabs}
\usepackage{hyperref}
\hypersetup{
    colorlinks=true,
    linkcolor=blue,
    filecolor=blue,      
    urlcolor=blue,
    citecolor=blue,
}

\theoremstyle{plain}
\newtheorem{thm}{Theorem}[section]

\newtheorem{pro}[thm]{Proposition}

\newtheorem*{claim*}{Claim}
\newtheorem{thmABC}{Theorem}

\theoremstyle{remark}

\newtheorem{dfn}[thm]{Definition}

\numberwithin{equation}{section}

\numberwithin{table}{section}

\newcommand{\Z}{\mathbb{Z}}
\newcommand{\Q}{\mathbb{Q}}
\newcommand{\F}{\mathbb{F}}

\newcommand{\mff}{\mathfrak{f}}

\newcommand{\mfp}{\mathfrak{p}}

\renewcommand{\epsilon}{\varepsilon}

\renewcommand{\phi}{\varphi}
\renewcommand{\theta}{\vartheta}

\newcommand{\mcO}{\mathcal{O}}

\newcommand{\Zp}{\mathbb{Z}_{p}}

\DeclareMathOperator{\gr}{gr}

\def \grL {\textup{gr}L}

\def \Fq {\ensuremath{\mathbb{F}_q}}

\DeclareMathOperator{\Gal}{Gal}

\theoremstyle{plain}
\newtheorem{mainthm}[thmABC]{Theorem}
\allowdisplaybreaks[2]
\hypersetup{
 pdftitle={Non-PORC ideal counting in two-generator free nilpotent Lie rings},
 pdfauthor={Marcus du Sautoy and Seungjai Lee},
 pdfsubject={Research announcement}
}

\ExplSyntaxOn
\NewCommandCopy \announcement_original_cite \cite
\RenewDocumentCommand \cite {o m}
 {
  \str_if_eq:nnTF {#2} {duS-ecI/01}
   {
    \IfNoValueTF {#1}
     {\announcement_original_cite{elliptic}}
     {\announcement_original_cite[#1]{elliptic}}
   }
   {
    \IfNoValueTF {#1}
     {\announcement_original_cite{#2}}
     {\announcement_original_cite[#1]{#2}}
   }
 }
\ExplSyntaxOff

\title[Non-PORC counting in free nilpotent Lie rings]{Non-PORC ideal counting in two-generator free nilpotent Lie rings}
\author{Marcus du Sautoy}\address{Mathematical Institute, University of Oxford, Oxford OX2 6GG, United Kingdom}\email{marcus.dusautoy@maths.ox.ac.uk}
\author{Seungjai Lee}\address{Department of Mathematics, Incheon National University, Incheon 22012, South Korea}\email{seungjai.lee@inu.ac.kr}
\date{26 November 2026}
\begin{document}
\begin{abstract}
We show that the normal subgroup zeta function of the free nilpotent group
$F_{8,2}$ of class eight on two generators is not a rational function on
residue classes of primes. More precisely, the number of its normal
subgroups of index $p^{13}$ is not polynomial on residue classes. This
gives a counterexample to the uniformity conjecture of Grunewald, Segal,
and Smith for normal subgroup zeta functions of free nilpotent groups.

We were led to this counterexample by counting graded ideals over finite
fields. In sufficiently large characteristic the number of graded ideals
of each codimension in $\mff_{c,2}(\Fq)$ is polynomial in $q$ through
class seven. In class eight the first exception occurs in codimension
thirteen, where the bracket conditions give the cubic $4z^3-10z^2-6z+9$.
The number of its roots modulo $p$ survives with nonzero coefficient
when we count all ideals of index $p^{13}$ in $\mff_{8,2}(\Zp)$.
The same root count occurs when the two-generator naturally graded Lie
algebras of dimension thirteen and class eight over $\F_p$ are counted
up to isomorphism. The Lazard correspondence thus gives a non-PORC
count for a family of groups of order $p^{13}$.
\end{abstract}
\maketitle

\section{Introduction}
\subsection{Zeta functions of free nilpotent groups and their uniformity}

Let $G$ be a finitely generated torsion-free nilpotent group
(a $\mathfrak{T}$-group for short). In their seminal paper \cite{GSS},
Grunewald, Segal and Smith introduced the \emph{normal subgroup zeta
function} of $G$,
\[
 \zeta_G^{\triangleleft}(s)
 =\sum_{m=1}^{\infty}a_m^{\triangleleft}(G)m^{-s}
 =\prod_p\zeta_{G,p}^{\triangleleft}(s),
 \qquad
 \zeta_{G,p}^{\triangleleft}(s)
 =\sum_{i=0}^{\infty}a_{p^i}^{\triangleleft}(G)p^{-is},
\]
where $a_m^{\triangleleft}(G)$ is the number of normal subgroups of
index $m$ in $G$. For each prime $p$, the local factor is a rational
function in $p^{-s}$ \cite[Theorem~1]{GSS}. The \emph{uniformity
problem} asks how this rational function varies with $p$.

\begin{dfn}
We say that $\zeta_G^{\triangleleft}(s)$ is \emph{finitely uniform}
if there exist $W_1,\ldots,W_k\in\Q(X,Y)$ such that, for each prime
$p$, the local factor is $W_i(p,p^{-s})$ for some $i$. It is
\emph{uniform} if a single rational function gives the local factors
for all sufficiently large primes, and \emph{non-uniform} if it is not
finitely uniform. Following \cite[Definition~1.1]{rforc}, we call it a
\emph{rational function on residue classes} (RFORC) if there exist
$M\geq1$ and $W_1,\ldots,W_M\in\Q(X,Y)$ such that, for all sufficiently
large primes $p$,
\[
 \zeta_{G,p}^{\triangleleft}(s)=W_i(p,p^{-s})
 \qquad\text{whenever }p\equiv i\pmod M.
\]
\end{dfn}

Free nilpotent groups occupied a special place in the original formulation
of the uniformity problem. Let $F_{c,d}$ denote the free nilpotent
group of class $c$ on $d$ generators. Grunewald, Segal and Smith
conjectured that a single rational function $W_{c,d}(X,Y)$ gives the
local normal subgroup zeta functions of $F_{c,d}$ for every prime
$p>c$ \cite[p.~188, Conjecture]{GSS}. They established uniformity in
classes one and two \cite[Proposition~1.1 and Theorem~2]{GSS}, and
Woodward obtained a uniform formula for $F_{3,2}$ \cite[Theorem~3.9]{Woodward/05}.

For a long time,  this conjecture seemed plausible. In particular, the first author asserted in \cite[Proposition~5.17(b)]{duSautoy2000}, and again in \cite[Proposition~3.19(b)]{duS-ennui/02}, that the conjecture holds for two generators in every class.

In this announcement article, we show that the conjecture is not true for $F_{8,2}$.
\begin{mainthm}\label{thm:non_uniform}
The normal subgroup zeta function $\zeta_{F_{8,2}}^{\triangleleft}(s)$
is not RFORC.
\end{mainthm}

The failure is already visible at index $p^{13}$. Recall that a function of
the prime $p$ is \emph{polynomial on residue classes} (PORC) if, for
some $M\geq1$, it agrees on each residue class modulo $M$ with a
polynomial in $p$ with rational coefficients. For a prime power $q=p^r$,
put
\begin{equation}\label{eq:f}
 f(z)=4z^3-10z^2-6z+9,
 \qquad N_q=\#\{z\in\F_q:f(z)=0\}.
\end{equation}

\begin{mainthm}\label{thm:p13_nonporc}
For all sufficiently large primes $p$,
\begin{equation}\label{eq:p13.Zp}
 a_{p^{13}}^{\triangleleft}(F_{8,2})
 =P_{\mathrm{ord}}(p)+p^{17}(p^2-1)(2p-1)N_p,
\end{equation}
where $P_{\mathrm{ord}}\in\Z[X]$ is given in~\eqref{eq:Pord}.
\end{mainthm}

Surprisingly, what led us to examine index $p^{13}$ was a much simpler counting
problem. We first encountered the cubic $f$ in the finite-field approximation
of Section~\ref{sec:approximations}. In sufficiently large characteristic
the graded ideal counts are polynomial in $q$ through class seven. In
class eight the first exception occurs in codimension thirteen, where
the count depends on $N_q$. Since a quotient of dimension thirteen over
$\F_p$ has order $p^{13}$, we looked for the same dependence on $N_p$
among ideals of this index over $\Zp$. Theorem~\ref{thm:p13_nonporc}
shows that the cubic term survives in the sum over all such ideals.
The variation of $N_p$ on residue classes (Section~\ref{sec:arithmetic})
then rules out uniformity for $F_{8,2}$.

\subsection{Zeta functions of free nilpotent Lie algebras and their approximations}\label{sec:approximations}

The finite-field approximation mentioned above combines two approaches to
counting ideals in nilpotent Lie rings. For a Lie algebra $L$ free of
finite rank over a commutative ring $R$, its ideal zeta function
$\zeta_L^{\triangleleft}(s)$ counts the $R$-ideals of finite index.
We write $a_m^{\triangleleft}(L)$ for the number of such ideals of
index $m$. By the Mal'cev correspondence, the local normal subgroup
zeta function of a $\mathfrak{T}$-group agrees, for all but finitely
many primes, with the local ideal zeta function of its associated Lie
ring \cite[Theorem~4.1]{GSS}. In particular, if $\mff_{c,d}(R)$
denotes the free nilpotent $R$-Lie algebra of class $c$ on $d$
generators, then
\[
 \zeta_{F_{c,d},p}^{\triangleleft}(s)
 =\zeta_{\mff_{c,d}(\Zp)}^{\triangleleft}(s)
\]
for all sufficiently large $p$.

Computing the local factors on the right is difficult in higher
nilpotency class. Voll and the second author \cite{gradedlocal}
studied an approximation that counts graded ideals. For a nilpotent
$R$-Lie algebra $L$ of class $c$, let
$L^{(i)}=\gamma_i(L)/\gamma_{i+1}(L)$ be the successive quotients
of its lower central series. The \emph{associated graded Lie algebra}
is $\grL=\bigoplus_{i=1}^c L^{(i)}$. An ideal $I$ of $\grL$ is
\emph{graded} or \emph{homogeneous} if
$I=\bigoplus_{i=1}^c(I\cap L^{(i)})$. Writing
$I\triangleleft_{\gr}\grL$ for a graded ideal of finite index, we
define the \emph{graded ideal zeta function} of $L$ by
\begin{equation}\label{def:graded.ideal.z.f.}
\zeta_{L}^{\triangleleft_{\gr}}(s)=\sum_{I\triangleleft_{\gr}\grL}|\grL:I|^{-s}.
\end{equation}
For $\mff_{c,d}(R)$ this is the grading by bracket length. Voll and
the second author obtained uniform graded ideal zeta functions for
$\mff_{3,3}(\Z)$ and $\mff_{4,2}(\Z)$
\cite[Theorem~1.1 and Proposition~5.3]{gradedlocal}.

A second approximation passes to the residue field. In an
eight-dimensional example studied by the second author, the number of
roots of $X^3-2$ modulo $p$ occurs both in the ideal zeta function
over $\F_p$ and in the local ideal zeta function over $\Zp$
\cite[Theorem~4.5]{Lee/25}.

We combine these approaches as follows. Let $\mcO_{\mfp}$ be the
completion of the ring of integers $\mcO$ of a number field at a
nonzero prime ideal $\mfp$, with residue field $\Fq=\mcO/\mfp$.
Put $L=\mff_{c,2}(\mcO_{\mfp})$, so that
$L/\mfp L\simeq\mff_{c,2}(\Fq)$. Taking inverse images under the
reduction map $L\to L/\mfp L$ gives a bijection between the graded
ideals of codimension $i$ in $L/\mfp L$ and the graded ideals of
index $q^i$ in $L$ containing $\mfp L$. We denote their number by
$a_{q^i}^{\triangleleft_{\gr}}(\mff_{c,2}(\Fq))$. For fixed $c$, we
call the graded ideal zeta function \emph{uniform} over finite fields
if, for each $i$, this number is given by a single polynomial in $q$
for all finite fields of sufficiently large characteristic.

\begin{mainthm}\label{thm:class_7}
For $1\leq c\leq7$, the graded ideal zeta function
$\zeta_{\mff_{c,2}(\Fq)}^{\triangleleft_{\gr}}(s)$ is uniform.
\end{mainthm}

Theorem~\ref{thm:class_7} corrects our earlier assertion that the
class-six coefficient in codimension nine depends on $p\bmod 8$
\cite[Theorems~1.6 and~6.1]{previous}. Correcting a bracket in degree
six gives a polynomial in $q$ for this coefficient in characteristic
greater than five.

\begin{mainthm}\label{thm:class_8}
Over finite fields of sufficiently large characteristic, the number of
graded ideals of $\mff_{8,2}(\Fq)$ in each codimension $0\leq i<13$
is polynomial in $q$. In codimension thirteen,
\[
 a_{q^{13}}^{\triangleleft_{\gr}}(\mff_{8,2}(\Fq))
 =P_8(q)+(q^3-q)N_q,
\]
where $P_8\in\Z[X]$ is given in~\eqref{eq:P8}. This coefficient is
not PORC, already for $q=p$ prime.
\end{mainthm}

To see how the same cubic enters the two counts, let $p$ be a
sufficiently large prime and put $L=\mff_{8,2}(\F_p)$.
Among the graded ideals $K$ for which $L/K$ has dimensions
$(2,1,2,2,2,2,1,1)$ in degrees one through eight, there is a family
whose bracket conditions reduce to a matrix in a parameter $v$ with
determinant $-(v+2)f(v)$. The roots of $f$ account for
$(p^3-p)N_p$ graded ideals in this family. Each such $K$ lifts to
$p^{16}(2p-1)$ ideals $I$ of index $p^{13}$ over $\Zp$. Here a lift
means that the images of the two free generators identify $L/K$ with
the Lie algebra associated to the lower central series of $\mff_{8,2}(\Zp)/I$.
These lifts therefore contribute
\[
 (p^3-p)N_p\,p^{16}(2p-1)=p^{17}(p^2-1)(2p-1)N_p .
\]
Summing the numbers of remaining ideals over all quotient classes gives
the polynomial $P_{\mathrm{ord}}(p)$, so this cubic term survives in
the full coefficient.

The polynomials in Theorems~\ref{thm:class_8} and~\ref{thm:p13_nonporc}
are
\begin{equation}\label{eq:P8}
\begin{aligned}
 P_8(X)={}&2X^{12}+6X^{11}+13X^{10}+27X^9+52X^8+87X^7\\
 &+124X^6+151X^5+169X^4+176X^3+139X^2+69X+26.
\end{aligned}
\end{equation}
\begin{equation}\label{eq:Pord}
\begin{aligned}
 P_{\mathrm{ord}}(X)={}&4X^{24}+26X^{23}+91X^{22}+257X^{21}+551X^{20}+815X^{19}+868X^{18}+729X^{17}\\&+565X^{16}+336X^{15}+167X^{14}+135X^{13}+119X^{12}+73X^{11}+36X^{10}+30X^9\\
 &+22X^8+15X^7+11X^6+7X^5+5X^4+3X^3+2X^2+X+1.
\end{aligned}
\end{equation}

\subsection{Quotients up to isomorphism and Higman's PORC conjecture}

In \cite{duSautoy2000} the first author established a strong connection between zeta functions of free nilpotent groups and the question of counting $p$-groups. Let $g_n(p)$ denote the number of isomorphism classes of groups of
order $p^n$. In 1960, Higman conjectured that $g_n(p)$ is PORC for each fixed
$n$ \cite{Higman/1960}, which is now known as Higman's PORC conjecture. Our graded ideal calculation gives the
number of isomorphism classes in one such family at order $p^{13}$.
Recall that a nilpotent Lie algebra $A$ is \emph{naturally graded} if
it is isomorphic to its associated graded Lie algebra for the lower
central series.

\begin{mainthm}\label{thm:graded-orbits}
For all sufficiently large primes $p$, the number of isomorphism
classes of two-generator naturally graded Lie algebras over $\F_p$
of dimension thirteen and nilpotency class eight is
\[
 p^3+6p^2+52p+175+2N_p .
\]
In particular, this number is not PORC.
\end{mainthm}

We obtain the formula by counting orbits of graded ideals under
linear changes of the two free generators. Each root of $f$ gives two
isomorphism classes in the family detected by the bracket conditions.
The two lines in the generator space determined by the degree-six
conditions are either both defined over $\F_p$ or interchanged by
Frobenius, giving one class in each case.
For $p>8$ the Lazard correspondence identifies these algebras with
the two-generator groups of order $p^{13}$, exponent $p$ and class
eight whose Lazard Lie algebra is naturally graded \cite{lazard}.
Theorem~\ref{thm:graded-orbits} therefore counts an explicitly
specified family in Higman's $g_{13}(p)$. 

\begin{mainthm}
The number of two-generated groups \(G\) of order \(p^{13}\), exponent \(p\), and class
at most eight for which \(\log(G)\) is naturally graded is not PORC.
\end{mainthm}

The other isomorphism
classes must also be counted to determine $g_{13}(p)$.

A full account of the proofs and computations announced here is in preparation.
\section{Splitting of the cubic and variation with the prime}\label{sec:arithmetic}
The reduction of $f$ modulo $5$ is $4(z^3+z+1)$ and has no root in
$\F_5$. Thus $f$ is irreducible over $\Q$. Its discriminant is
$46944=2^5\cdot3^2\cdot163$, so its splitting field $K$ has Galois
group $S_3$. At each prime of good reduction, $N_p$ is the number
of fixed points of Frobenius on the three roots.
Since $K/\Q$ is nonabelian, the criterion for complete splitting
\cite[Theorem~2.2.1]{Weinstein/2016} shows that the condition $N_p=3$
cannot be determined by the residue class of $p$ modulo any fixed integer.

\begin{pro}\label{prop:nonporc}
For every $m\geq1$, each of $N_p=0$ and $N_p=3$ occurs for infinitely
many primes $p\equiv1\bmod m$. If $A,B\in\Q(X)$ and $B\ne0$, then
$A(p)+B(p)N_p$ cannot agree with a single rational function of $p$
on all sufficiently large primes in this progression.
\end{pro}
\begin{proof}
Put $C=\Q(\zeta_m)$. Since $K\cap C$ is abelian over $\Q$, the
commutator subgroup $A_3$ of $\Gal(K/\Q)$ fixes $K\cap C$.
Thus the identity and a three-cycle on $K$ each extend to an
automorphism of $KC$ acting trivially on $C$. By Chebotarev's theorem
\cite[Chapter~VII, Theorem~13.4]{Neukirch/99}, each of the two conjugacy classes in
$\Gal(KC/\Q)$ so obtained occurs as Frobenius at infinitely many primes.
These primes satisfy $p\equiv1\bmod m$ because they split completely
in $C$. The identity fixes all three roots of $f$, and a three-cycle
fixes none. A rational function $R$ as in the statement would therefore
agree with $A+3B$ at infinitely many primes and with $A$ at infinitely
many others, giving $R=A+3B=A$ in $\Q(X)$, contrary to $B\ne0$.
\end{proof}

At $q=p$, the proposition proves the non-PORC assertions in
Theorems~\ref{thm:p13_nonporc}, \ref{thm:class_8}
and~\ref{thm:graded-orbits}.
If the local factors were RFORC with modulus $M$, expanding the
rational function for $p\equiv1\pmod M$ in $T=p^{-s}$ would give
the coefficient of $T^{13}$ as a rational function of $p$ for all
sufficiently large primes in that progression.
Proposition~\ref{prop:nonporc} rules this out for both coefficients
in Theorems~\ref{thm:p13_nonporc} and~\ref{thm:class_8}.

\subsection*{AI Declaration}
The authors used ChatGPT 6.2 Astra in developing mathematical arguments, computing  explicit calculations, and running verification programs.

\bibliographystyle{amsplain}
\bibliography{references}
\end{document}